\documentclass[11pt]{amsart}
\usepackage[vmargin=1in,hmargin=1.25in]{geometry}                
\usepackage{graphicx}
\usepackage{amssymb}
\usepackage{epstopdf}
\usepackage{url}
\usepackage{amsmath,amsthm}
\title{Gauss and the Eccentric Halsted}
\author{Jonathan Sondow}

\begin{document}
\begin{abstract}
We discuss the mathematician George Bruce Halsted's accusations against Carl Friedrich Gauss, as well as refutations both by the latter's American grandson Robert Gauss in a letter to Felix Klein, and by the historian of mathematics Florian Cajori.
\end{abstract}
\maketitle
\thispagestyle{empty}




In the January, 1912, issue of the \emph{Monthly} \cite{halsted}, the mathematician and translator George Bruce Halsted (1853--1922) had this to say about Carl Friedrich Gauss (1777--1855):
\begin{quote}
\ \quad It has been said the heroic age of non-euclidean geometry is passed, since long gone are the days when Lobachevski flinched into calling his system ``imaginary geometry,'' and, I might add, Gauss kept a more cowardly silence because, as he confesses, he ``fears the outcry of the Boeotians,'' (proverbial for stupidity).

\dotso the designation ``Absolute Geometry'' was first used by myself as a rendering for John Bolyai's phrase \emph{scientiam spatii absolute veram}, in which I have always gloried as showing the magnificent verve of the young Magyar hero, victim of the meanness of Gauss, as was also his own son who passed his life an exile here in Colorado.
\end{quote}

\centerline{\includegraphics[width=3.7in]{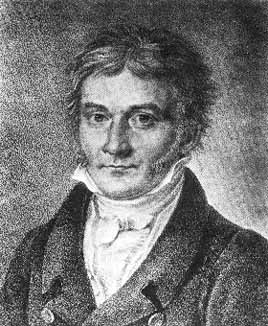}}\begin{center}{Carl Friedrich GAUSS (1777--1855)}\end{center}

\centerline{\includegraphics[width=3in]{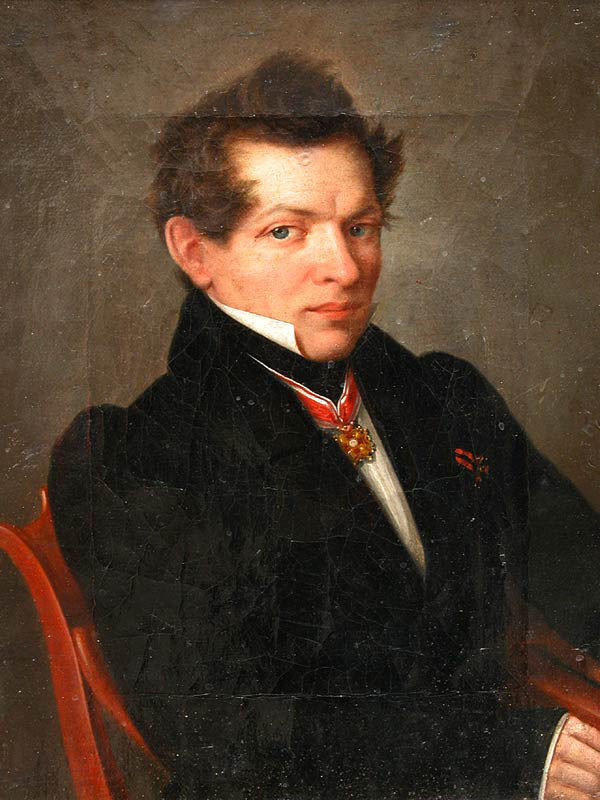}}\begin{center}{Nikolai Ivanovich LOBACHEVSKY (1853--1922)}\end{center}

\centerline{\includegraphics[width=2.8in]{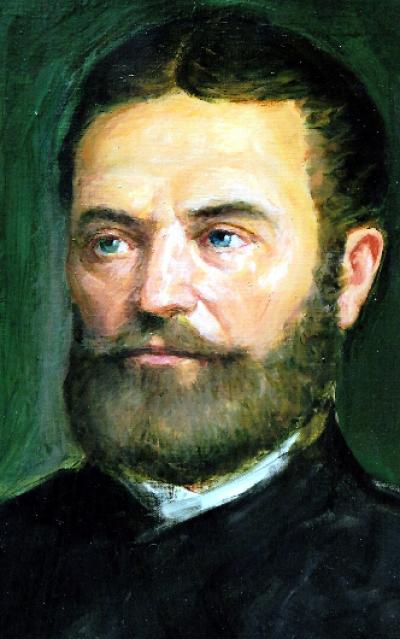}}\begin{center}{J\'{a}nos BOLYAI (1802--1860)}\end{center}

\centerline{\includegraphics[width=3in]{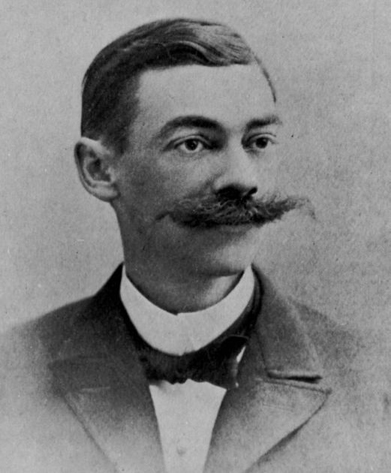}}\begin{center}{George Bruce HALSTED (1792--1856)}\end{center}

In April, 1912, Halsted issued the following ``correction'' \cite{halsted12}:
\begin{quote}
\ \quad My friend, Professor Cajori, informs me that no one of Gauss's four sons made his home in Colorado. He tells me that the son Eugen had a quarrel with his father, because of Eugen's wild life at the University of G\"{o}ttingen, but when Eugen left home to sail for America, his father followed him, urged him to return home, and, when failing in this effort, offered him money to take with him to America. Eugen left home and came to this country by his own choice and against the wishes of his parents. His father was greatly grieved, as is shown by his correspondence. See the article \dotso \cite{cajori99}.
\end{quote}

\centerline{\includegraphics[width=2.6in]{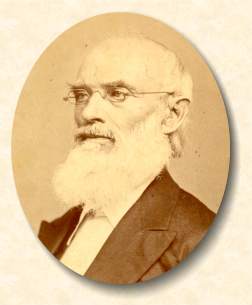}}\begin{center}{Eugene GAUSS (1811--1896)}\end{center}

In September, 1912, Halsted's personal accusations against Gauss were refuted  by the latter's American grandson Robert Gauss (1851--1913) in a letter \cite{gauss} to Felix Klein (1849--1925). (Klein's championing of Gauss's role in the development of non-Euclidean geometry---``There can be
no doubt that Lobach\'{e}vski as well as Bolyai owe to Gauss's prompting the initiative
of their researches.''---had been ridiculed by Halsted, who predicted, ``The whole scientific world will breath a sigh of relief that Klein's Goettingen legend, wounded in 1893, is in 1899 annihilated forever.'' \cite{halsted99}) In his letter to Klein, Robert Gauss wrote:

\begin{quote}
\quad \quad \dotso I am sure that you know by reputation, if not personally, Mr. George Bruce Halsted, who has given much attention to the study of Non-Euclidean Geometry, a branch of mathematics in which you have taken great interest. If I am correct in this, you are aware that he has been a great champion of John Bolyai. For some reason this has developed in him a spirit of antagonism to my grandfather and led him into a very unjust attack upon him. In the January number for this year of ``The American Mathematical Monthly'' he had an article in which he declared that John Bolyai was a victim of ``the meanness of Gauss''. To this he added that one of Gauss' own sons also was a victim of this ``meanness'', and that he had spent his life an exile in the State of Colorado.

Professor Florian Cajori of Colorado Springs in this State, wrote me after reading Professor Halsted's article. He knew that the reference could only be to Eugene Gauss, my father. \dotso I promptly replied to Professor Cajori that my father was in no sense an exile from his home in G\"{o}ttingen, and furthermore, that he had never been in the State of Colorado.

\dotso The enclosed letters also show the circumstances under which my father left G\"{o}ttingen, and I think they are a complete vindication of my grandfather against the charge made by Halsted that he had ``meanly'' treated my father.

\dotso I should have said that Professor Cajori wrote an answer \cite{cajori12} to Halsted's attack upon my grandfather in relation to Bolyai and conclusively refuted it.  In that article he also made some reference to my father \dotso.
\end{quote}

\centerline{\includegraphics[width=2.5in]{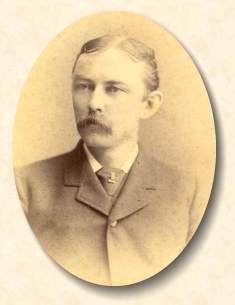}}\begin{center}{Robert GAUSS (1851--1913)}\end{center}

\centerline{\includegraphics[width=2.6in]{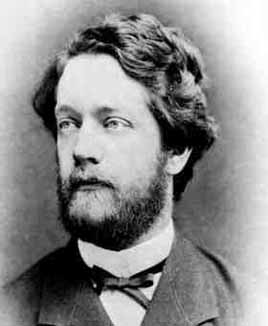}}\begin{center}{Felix Christian KLEIN (1849--1925)}\end{center}

A decade later, in an obituary \cite{cajori} for Halsted, the historian of mathematics Florian Cajori (1859--1930) summed up the affair this way:

\begin{quote}
\ \quad American mathematicians are indebted to Halsted for making the writings of the creators of non-Euclidean geometry accessible to them in the English language. His commentaries were always spicy and valuable, even though, as a historian, Halsted was not always able to maintain the attitude of an impartial judge. At his hands Gauss, for instance, received scant justice.
\end{quote}

Further writings by Cajori and Halsted on Gauss and non-Euclidean geometry can be found in \cite{cajori12, halsted00, halsted04}. Revealing biographies of ``the eccentric and sometimes spectacular Halsted'' are in \cite{or, tropp}, an anecdotal one of the ``erratic, and somewhat bizarre Halsted'' is in \cite{benedict}, and an appreciation of ``the inexhaustible fruitfulness of Dr. Halsted's genius'' by his student L.~E. Dickson is in \cite{dickson}; see also the section of \cite{zitarelli} on Halsted's student R. L. Moore.

\centerline{\includegraphics[width=2.1in]{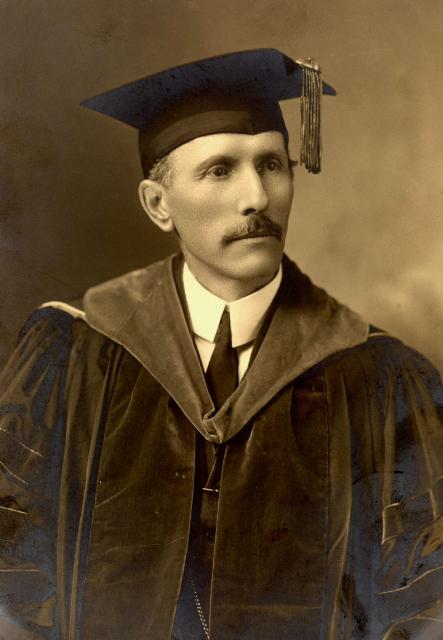}}\begin{center}{Florian CAJORI (1859--1930)}\end{center}

\smallskip
\noindent\textit{209 West 97th Street, New York, NY 10025\\
jsondow@alumni.princeton.edu}

\end{document}